\documentclass[11pt, reqno]{amsart} 
\usepackage{amssymb, latexsym, amsmath}
\newtheorem{theorem}{Theorem}[section]
\newtheorem{lemma}[theorem]{Lemma}
\newtheorem{proposition}[theorem]{Proposition}

\newtheorem{corollary}[theorem]{Corollary}
\theoremstyle{definition}
\newtheorem{remark}[theorem]{Remark} 
\theoremstyle{definition}
\newtheorem{definition}[theorem]{Definition} 
\theoremstyle{definition}
\newtheorem{example}[theorem]{Example}
\theoremstyle{definition}
\newtheorem{construction}[theorem]{Construction}
 
\numberwithin{equation}{section}

\usepackage[mathscr]{eucal}

\usepackage{xy} 
\xyoption{all}

\newcommand{\PP}{\mathbb{P}}
\newcommand{\wPP}{\widetilde{\mathbb{P}}}
\newcommand{\CC}{\mathbb{C}}
\newcommand{\ZZ}{\mathbb{Z}}

\newcommand{\RR}{\mathbb{R}}
\newcommand{\TT}{\mathbb{T}}
\newcommand{\mc}[1]{\mathcal{#1}}

\newcommand{\mb}[1]{\mathbb{#1}}
\newcommand{\del}{\partial}
\newcommand{\Oh}{\mathcal{O}}
\newcommand{\E}{\mathcal{E}}

\DeclareMathOperator{\Spec}{Spec}
\DeclareMathOperator{\Sym}{Sym}

\DeclareMathOperator{\sing}{sing}
\DeclareMathOperator{\Hom}{Hom}
\DeclareMathOperator{\Bs}{Bs}

\DeclareMathOperator{\Proj}{Proj}

\DeclareMathOperator{\tr}{tr}
\DeclareMathOperator*{\Span}{Span}

\DeclareMathOperator{\sheafHom}{\mathscr{H}\mspace{-1mu}\mathit{om}}

\DeclareMathOperator{\codim}{codim}

\CompileMatrices

\begin{document}
\title[Toric Calabi-Yau Hypersurfaces]{On toric Calabi-Yau hypersurfaces
fibered by weighted K3 hypersurfaces}
\date{\today}

\author[J. Mullet]{Joshua P. Mullet}
\address{Department of Mathematics\\
        The Ohio State University\\
        Columbus, OH 43210}
\email{mullet@math.ohio-state.edu}

\begin{abstract}
In response to a question of Reid, we find all anti-canonical Calabi-Yau
hypersurfaces $X$ in toric weighted $\PP^3$-bundles over $\PP^1$ where
the general fiber of $X$ over $\PP^1$ is a weighted K3 hypersurface.
This gives a direct generalization of Reid's discovery of the 95
families of weighted K3 hypersurfaces in \cite{mR80}.  We also treat the
case where $X$ is fibered over $\PP^2$ with general fiber a 
genus one curve in a weighted projective plane.
\end{abstract}
\maketitle

\section{Introduction} 
In light of the successes of the minimal model program in the
classification of threefolds, Corti and Reid promulgate in \cite{aC00} a
program of \emph{explicit} birational study of threefolds.  Simply put,
the interest is to obtain explicit equations for threefolds as
hypersurfaces or complete intersections in simple varieties such as
projective spaces or weighted projective spaces, and then focus study to
these explicit examples, regarding them as birational models.  An early
example of such an explicit result is Reid's discovery of the 95
families of K3 hypersurfaces in weighted projective 3-space (see
\cite{mR80}).  

Miles Reid posed to the author the problem of classifying Calabi-Yau
threefolds fibered over the projective line where the general fiber is
one of the 95 weighted K3 hypersurfaces mentioned in the previous
paragraph. Since the fibers are to be weighted hypersurfaces, a natural
place to look for such threefolds is in weighted projective bundles.
Furthermore, if we are interested in being explicit, we should look for
these threefolds in some ``simple'' varieties.  For us, those simple
varieties are toric varieties.  We are thus naturally led to the problem
of finding K3-fibered Calabi-Yau threefolds in toric weighted projective
bundles over $\PP^1$.  In this paper, we show how to find all such
threefolds and thereby obtain a direct generalization of Reid's
discovery of the 95 families of weighted K3 hypersurfaces.

The contents of this paper are as follows. In
Section~\ref{S:conventions} we fix some definitions and notation
regarding the theory of toric varieties that we use throughout the
paper.  In Section~\ref{S:wpsbund}, we construct weighted projective
space bundles as toric varieties.  Weighted projective space bundles can
also be obtained by taking $\Proj$ of a sheaf of weighted polynomial
algebras.  In Section~\ref{S:proj} we describe the $\Proj$ construction
and compare it to the construction using toric varieties.  In
Section~\ref{S:intth} we describe the intersection theory of weighted
$\PP^1$-bundles over $\PP^1$ in terms of the well-known intersection
theory on ordinary $\PP^1$-bundles.  We will need this description in
the proof of Theorem~\ref{T:qsmoothk3fib}.  In Section~\ref{S:linsys},
we consider linear systems on toric varieties and, more specifically,
on weighted projective spaces and bundles.  Section~\ref{S:wform} is a
technical section in which we show that the condition of well-formedness
(Definition~\ref{D:wformsub}), which is a condition imposed on
hypersurfaces to make the adjunction formula work as expected, is
automatic in the case of quasi-smooth hypersurfaces
(Definition~\ref{D:quasismooth}).  In Section~\ref{S:achs}, we prove our
main result, Theorem~\ref{T:qsmoothk3fib}, in which we give necessary
and sufficient conditions for a toric weighted projective bundle over
$\PP^1$ to admit a K3-fibered anti-canonical Calabi-Yau hypersurface,
thus providing an answer to Reid's original question.
Section~\ref{S:ellfib} contains a statement of the analogous result for
elliptically fibered Calabi-Yau threefolds over $\PP^2$ whose general
fiber is a genus one curve in a weighted $\PP^3$.  The complete list of
all our Calabi-Yau varieties is quite long, but the appendix contains a
sampling of data, and we refer the reader to \cite{jM06} for all the
data and details of the calculation. 

We should mention that Masanori Kobayashi considered the problem in the
non-weighted case (\cite{mK}) and without the explicit use of toric
varieties.  Our work is inspired by his, and we use some of his notation
in Section~\ref{S:linsys}.

\subsection*{Acknowledgments}\label{Ss:acknowledgments}
I would like to thank my advisers, Dan Grayson and Sheldon Katz, for
their attention and encouragement.  I would also like to
thank Miles Reid for suggesting the problem and Masanori Kobayashi
for sharing unpublished work with me. 
    
\section{Conventions and notation}\label{S:conventions}
We work over $\CC$, the field of complex numbers, and we
follow standard conventions regarding the theory of toric varieties as
described, for example, in \cite{wF93} or \cite{vD78}.  If $X$ is a
toric variety, we will denote its associated lattices as $N_X$ and $M_X
:= \Hom_{\ZZ}(N_X, \ZZ)$, and if $N$ is a lattice we define $N_\RR := N
\otimes_\ZZ \RR$.  We always assume that our toric varieties have the
property that all their maximal cones have dimension equal to the rank
of the lattice used to define them.

We will need the notions of \emph{quasi-smoothness}
and \emph{well-formedness}, which we define below for subvarieties of
toric varieties.  Let $X_\Sigma$ be a simplicial toric variety
associated to a fan $\Sigma$, and let $S_\Sigma$ be its homogeneous
coordinate ring with irrelevant ideal $B \subseteq S_\Sigma$ (see
\cite{dC95}).  According to \cite[Theorem~2.1]{dC95}, $X_\Sigma$ is
isomorphic to the geometric quotient
\[
    (\Spec S_\Sigma - Z(B))/G
\]
where 
\[
    G = \Hom_\ZZ(A^1X_\Sigma, \CC^*).
\]

\begin{definition}\label{D:quasismooth}
Let 
\[
\xymatrix{
    {\Spec S_\Sigma - Z(B)} \ar[r]^-*{q} & {X_\Sigma} }
\]
be the quotient map.  We say that a closed subvariety $Y \subseteq
X_\Sigma$ is \emph{quasi-smooth} if $q^{-1}Y$ is nonsingular.
\end{definition}

\begin{definition}\label{D:wformsub}
Let $X_\Sigma$ be a simplicial toric variety, and let $\sing X_\Sigma$
denote the singular locus of $X_\Sigma$.  We say that a closed
subvariety $Y \subseteq X_\Sigma$ is \emph{well-formed} if the
codimension of $(\sing X_\Sigma) \cap Y$ is at least two in $Y$.
\end{definition}

\begin{remark}\label{R:orbifold}
If $Y \subseteq X_\Sigma$ is quasi-smooth, then $Y$ is an orbifold.
\end{remark}

\section{Weighted Projective Space Bundles as Toric
Varieties}\label{S:wpsbund}

In this section we consider the problem of constructing weighted
projective bundles as toric varieties.  We first recall how to
construct the weighted projective space 
\[
    \PP(a_0, \ldots, a_n)
\]
as a toric variety. 

\begin{construction}\label{Co:wpstv}
Let $(a_0, \ldots, a_n)$ be a list of
positive integers such that 
\begin{equation}\label{Eq:wform}
    \gcd(a_0, \ldots, \widehat{a_i}, \ldots, a_n) = 1
\end{equation}
where $\widehat{a_i}$ means that we omit $a_i$.  Then form the weighted
projective space 
\[
    \PP(a_0, \ldots, a_n) := \Proj \CC[X_0, \ldots, X_n]
\]
where $\deg X_i = a_i$. (See \cite{iD82, aI00} for generalities about
weighted projective space.)  As described in \cite{wF93}, we may write
$\PP(a_0, \ldots, a_n)$ as a toric variety as follows.  Define the
lattice $N$ via the following exact sequence:
\begin{equation}\label{Eq:wpses}
\xymatrix{
{0} \ar[r] & {\ZZ} \ar[r]_-*{
        \left(
        \begin{smallmatrix}
        a_0 \\ \vdots \\ a_n
        \end{smallmatrix}
        \right)}
            & {\ZZ^{n+1}} \ar[r] & {N} \ar[r] & {0},}
\end{equation}
and define the fan $\Sigma \subseteq N_\RR$ to consist of
cones spanned by the images of proper subsets of the standard basis for
$\ZZ^{n+1}$.  Then
\[
    \PP(a_0, \ldots, a_n) \cong X_\Sigma.
\]
Note that if we take $a_i = 1$ for all $0 \leq i \leq n$, we obtain a
description of (ordinary) projective space as a toric variety.
\end{construction}

We come now to the notion of fibration for toric varieties that we will
use (cf. \cite[Exercise, p.~41]{wF93}). 

\begin{definition}\label{D:toricfib}
By a \emph{locally trivial toric fibration} we mean a fibration
\begin{equation*}\label{Eq:toricfib}
\xymatrix{
    {F} \ar[r]^*{i} & {X} \ar[d]_*{\pi} \\
    {} & {B}  \\ }
\end{equation*}
in which $F$, $X$, and $B$ are toric varieties associated to fans
$\Sigma_F$, $\Sigma_X$, $\Sigma_B$ and lattices $N_F$, $N_X$, $N_B$
respectively.  We further require that the fibration trivialize over
affine open toric subvarieties  of $B$, that $F$ is the fiber over the
identity element of the torus for $B$, and that the maps $i$ and $\pi$ 
are induced by maps of fans $i'$ and $\pi'$ respectively.
\end{definition} 

\begin{remark}\label{R:vulgar}We may be imprecise and say that $X$ is a
``fibration'' with fiber $F$, or a ``toric $F$-bundle.'' 
\end{remark}

\begin{remark}\label{R:bes}
Note that the maps of lattices in Definition~\ref{D:toricfib} form  an
exact sequence
\[
\xymatrix{
    {0} \ar[r] & {N_F} \ar[r]^*{i'} & {N_X} \ar[r]^*{\pi'} & {N_B}
    \ar[r] & {0}. }
\] 
\end{remark}

\begin{construction}\label{Co:toricwpb}
Let $n$ and $k$ be positive integers, let $(a_0, \ldots, a_n)$ be a
list of positive integers satisfying \eqref{Eq:wform}, and let
$(d_0, \ldots, d_n)$ be an arbitrary list of integers. Define a lattice
$N$ via the following exact sequence:
\begin{equation}\label{Eq:toricwpb}
\xymatrix{
    {0} \ar[r] & {\ZZ^2} \ar[rr]_-*{
        \left(
        \begin{smallmatrix}
            a_0 & -d_0 \\
            \vdots & \vdots \\
            a_n & -d_n \\
            0 & 1 \\
            \vdots & \vdots \\
            0 & 1
        \end{smallmatrix}
        \right)
    } && {\ZZ^{n+1} \oplus \ZZ^{k+1}} \ar[r] & {N} \ar[r] & {0} 
}.
\end{equation}
Now let $\{e_0, \ldots, e_n\}$ denote the standard basis for
$\ZZ^{n+1}$, and let $\{f_0,\ldots, f_k\}$ denote the standard basis for
$\ZZ^{k+1}$.  We let $\Sigma \subseteq N_\RR$ be the fan that consists
of cones spanned by proper subsets of the standard bases that contain
neither all of the $e_i$ nor all of the $f_i$.  
\end{construction}

\begin{proposition}\label{P:istoricbund}
Let $\Sigma$ be the fan constructed in Construction~\ref{Co:toricwpb},
and let $X_\Sigma$ be the associated toric variety.  Then $X_\Sigma$ is
a locally trivial toric fibration over $\PP^k$ (in the sense of
Definition~\ref{D:toricfib}) over $\PP^k$ with fiber $\PP(a_0, \ldots,
a_n)$.
\end{proposition}

\begin{proof}
Let $N_{X_\Sigma}$ be the lattice for $X_\Sigma$, let $N_{\PP^k}$ be
the lattice for $\PP^k$, and let $N_{\PP(a_0,\ldots,a_n)}$ be the
lattice for $\PP(a_0, \ldots, a_n)$; here we regard $\PP^k$ and
$\PP(a_0, \ldots, a_n)$ as toric varieties according to
Construction~\ref{Co:wpstv}.  We have the following commutative diagram
whose rows and columns are exact.
\begin{equation}\label{Eq:istoricdia}
\xymatrix{
    {} & {0} \ar[d] && {0} \ar[d] & {0}  \ar[d] & {} \\
    {0} \ar[r] & {\ZZ} \ar[d] \ar[rr]_-*{
        \left(
        \begin{smallmatrix}
        a_0 \\ \vdots \\ a_n
        \end{smallmatrix}
        \right)
    } && {\ZZ^{n+1}} \ar[d] \ar[r] & {N_{\PP(a_0, \ldots, a_n)}} \ar[d]
    \ar[r] & {0} \\ 
    {0} \ar[r] & {\ZZ^2} \ar[dd] \ar[rr]_-*{
        \left(
        \begin{smallmatrix}
            a_0 & -d_0 \\
            \vdots & \vdots \\
            a_n & -d_n \\
            0 & 1 \\
            \vdots & \vdots \\
            0 & 1
        \end{smallmatrix}
        \right)
    } && {\ZZ^{n+1} \oplus \ZZ^{k+1}} \ar[dd] \ar[r] & {N} \ar[dd]
    \ar[r] & {0} \\ \\
    {0} \ar[r] & {\ZZ} \ar[d] \ar[rr]_-*{
        \left(
        \begin{smallmatrix}
        1 \\ \vdots \\ 1
        \end{smallmatrix}
        \right)
    } && {\ZZ^{k+1}} \ar[d] \ar[r] & {N_{\PP^k}} \ar[d] \ar[r] & {0} \\
    {} & {0} && {0} & {0} & {} 
}
\end{equation}
The right column of \eqref{Eq:istoricdia} gives the required maps of
lattices.  Furthermore, since the middle column is just inclusion of the
first factor followed by projection onto the second factor, we can
canonically extend the maps of lattices to maps of fans.  These maps of
fans give the required maps of toric varieties.
\begin{equation}\label{Eq:wpbfib}
\xymatrix{
    {\PP(a_0, \ldots, a_n)} \ar[r]^-*{i} & {X_\Sigma} \ar[d]_*{\pi} \\
    {} & {\PP^k}  
}
\end{equation}

It remains to show that \eqref{Eq:wpbfib} trivializes over affine open
toric subvarieties of $\PP^k$, and for this it suffices to
consider only those open toric subvarieties of $\PP^k$ that correspond
to maximal cones in the fan for $\PP^k$.  Let $U_i$ be the open toric
subvariety corresponding to the cone spanned by the images of all but
the $i$-th standard basis vector of $\ZZ^{k+1}$.  Next observe that for
all $0 \leq i \leq k$ we have an isomorphism of exact sequences 
\[
\xymatrix{
    {0} \ar[r] & {\ZZ^2} \ar[dd] \ar[rr]^-*{
        \left(
        \begin{smallmatrix}
            a_0 & -d_0 \\
            \vdots & \vdots \\
            a_n & -d_n \\
            0 & 1 \\
            \vdots & \vdots \\
            0 & 1
        \end{smallmatrix}
        \right)
    } && {\ZZ^{n+1} \oplus \ZZ^{k+1}} \ar[dd]_-*{
        \left(
        \begin{smallmatrix}
            I_{n+1} & M \\
            0 & I_{k+1}
        \end{smallmatrix}
        \right)
    } \ar[r] & {N} \ar[dd]_-*{\cong} \ar[r] & {0} \\ \\
    {0} \ar[r] & {\ZZ^2} \ar[rr]_-*{
        \left(
        \begin{smallmatrix}
            a_0 & 0 \\
            \vdots & \vdots \\
            a_n & 0 \\
            0 & 1 \\
            \vdots & \vdots \\
            0 & 1
        \end{smallmatrix}
        \right)
    } && {\ZZ^{n+1} \oplus \ZZ^{k+1}} \ar[r] & {N'} \ar[r] & {0}  
}
\]
where $M$ is an $(n+1) \times (k+1)$ matrix whose only nonzero column
lies in the same column as the $i$-th basis vector in $I_{k+1}$ and
is equal to
\[
    \left(
    \begin{smallmatrix}
    d_0 \\ \vdots \\ d_{n}
    \end{smallmatrix}
    \right).
\]
Therefore, if we consider the open subvariety of $X_\Sigma$ that
corresponds to cones that do not contain the ray corresponding to the
$i$-th standard basis vector of $\ZZ^{k+1}$, the map of exact sequences
above yields an isomorphism of fans for $\pi^{-1}U_i$ and $\PP(a_0,
\ldots, a_n) \times U_i$ where $\pi$ is the structure map in
\eqref{Eq:wpbfib}.
\end{proof}

Our next result shows that all toric weighted projective space
bundles are of the form described in Construction~\ref{Co:toricwpb}.  

\begin{theorem}\label{T:wpsbund/ps}
Let 
\begin{equation}\label{Eq:wpsbund/ps}
\xymatrix{
    {\PP(a_0, \ldots, a_n)} \ar[r]^-*{i} & {\wPP} \ar[d]_*{\pi} \\
    {} & {\PP^k} }
\end{equation}
be a locally trivial toric fibration as in Definition~\ref{D:toricfib},
where we regard $\PP^k$ as a toric variety according to
Construction~\ref{Co:wpstv}.  Then there exist integers $d_0, \ldots,
d_n$ such that \eqref{Eq:wpsbund/ps} is isomorphic to the fibration
constructed in Construction~\ref{Co:toricwpb}. 
\end{theorem}

\begin{proof}
Let $N_{\wPP}$ be the lattice for ${\wPP}$.  Let $\rho_0, \ldots,
\rho_n$ be the rays of the fan for $\PP(a_0, \ldots, a_n)$, and let
$\tau_0, \ldots, \tau_k$ be the rays of the fan for $\PP^k$.  Following
standard convention, we write $v_{\rho_i}$ and $v_{\tau_i}$ for the
primitive vectors in the rays $\rho_i$ and $\tau_i$.  Since the
fibration is locally trivial, there exist unique rays
$\widetilde{\tau_0}, \ldots, \widetilde{\tau_k}$ in $N_{\wPP}$ such that
$\pi'(v_{\widetilde{\tau_i}}) = v_{\tau_i}$ for $0 \leq i \leq k$
where $\pi'$ is the map of fans corresponding to the map $\pi$.  Let
$\{e_0, \ldots, e_n\}$ be the standard basis for $\ZZ^{n+1}$ and let
$\{f_0, \ldots, f_k\}$ be the standard basis for $\ZZ^{k+1}$.  We define a
map 
\[
\xymatrix{
    {F:\ZZ^{n+1} \oplus \ZZ^{k+1}} \ar[r] & {N_{\wPP}} }
\]
by setting
\[
    F(e_i,f_j) := i'(v_{\rho_i}) + v_{\widetilde{\tau_j}}
\]
where $i'$ is the map of fans corresponding to the map $i$ in
\eqref{Eq:wpsbund/ps}.  We easily check that $F$ is surjective, and we
have an exact sequence,
\begin{equation}\label{Eq:prodes}
\xymatrix{
    {0} \ar[r] & {\ZZ^2} \ar[r] & 
        {\ZZ^{n+1} \oplus \ZZ^{k+1}} \ar[r]^-*{F} & {N_{\wPP}} \ar[r] &
        {0}}.
\end{equation}
Since the fibration is locally trivial, we may use the fact
(\cite[Exercise, p. 22]{wF93}) that the fan for the product of two toric
varieties is isomorphic to the product of their fans to check that the
fan $\Sigma_{\wPP}$ is as described in Construction~\ref{Co:toricwpb}
but now with respect to the map $F$.  

It remains to show that the kernel of \eqref{Eq:prodes} is isomorphic to
the kernel in \eqref{Eq:toricwpb}. From \eqref{Eq:wpses} we see that
$\sum a_i F(e_i) = 0$.  Furthermore, since $\sum v_{\tau_i} = 0$,  we
can apply Remark~\ref{R:bes} to find
\[
    F\left(\sum f_i\right) = F\left(\sum d_j e_j\right)
\] 
for some integers $d_0, \ldots, d_n$.  Therefore, \eqref{Eq:prodes} fits
into a commutative diagram
\[
\xymatrix{
    {} & {\ZZ^2} \ar@{^{(}->}[d]_*{\phi} \ar[dr]^*{
        \left(
        \begin{smallmatrix}
            a_0 & -d_0 \\
            \vdots & \vdots \\
            a_n & -d_n \\
            0 & 1 \\
            \vdots & \vdots \\
            0 & 1 
        \end{smallmatrix}
        \right) } &&& \\
    {0} \ar[r] & {\ZZ^2} \ar[r] & 
        {\ZZ^{n+1} \oplus \ZZ^{k+1}} \ar[r] & {N_{\wPP}} \ar[r] & {0}\ ,}
\]
but by \eqref{Eq:wform} the columns in the matrix are primitive,
therefore $\phi$ is an isomorphism as required.
\end{proof}

Having characterized all toric $\PP(a_0, \ldots, a_n)$-bundles over
$\PP^k$, we conclude this section by finding their homogeneous
coordinate rings and irrelevant ideals (\cite{dC95}).  The description
given in Construction~\ref{Co:toricwpb} makes this particularly easy.
Recall that for a toric variety $X$ there exists an exact sequence
\begin{equation}\label{Eq:torices}
\xymatrix{
    {0} \ar[r] & {M} \ar[r] & {\ZZ^{\Sigma(1)}} \ar[r] & {A^1} X \ar[r] &
    {0,}}
\end{equation}
where $\Sigma(1)$ denotes the set of one-dimensional cones in the fan,
and where $A^1 X$ denotes the Chow group of codimension one cycles
modulo rational equivalence.
We begin
with a proposition in which we compute the exact sequence
\eqref{Eq:torices} for toric weighted projective bundles.

\begin{proposition}\label{P:wpstes}
Let $(a_0, \ldots, a_n)$ be a sequence of positive integers
satisfying \eqref{Eq:wform}, and let ${\wPP}$ be a locally trivial
toric fibration over $\PP^k$ with fiber $\PP(a_0, \ldots, a_n)$.  There
is an isomorphism of exact sequences:
\begin{equation*}\label{Eq:wpstes}
\xymatrix{
    {0} \ar[r] & {M_{\wPP}} \ar[r] \ar@{=}[d] & {\Hom(\ZZ^{n+1} \oplus
    \ZZ^{k+1}, \ZZ)} \ar[d] \ar[r]^-*{T^{\tr}} & {\Hom(\ZZ^2, \ZZ)}
    \ar[r] \ar[d] & {0} \\
    {0} \ar[r] & {M_{\wPP}} \ar[r] & {\ZZ^{\Sigma(1)}} \ar[r] &
    {A^1({\wPP})}
        \ar[r] & {0}\ , }
\end{equation*} 
where $T$ is the matrix in Construction~\ref{Co:toricwpb} and the bottom
row is the standard exact sequence \eqref{Eq:torices}.
\end{proposition}

\begin{proof}
Since the rays of the fan for ${\wPP}$ are given by the images of
standard basis vectors under the cokernel of the map $T$, the result
will follow immediately provided we show that the images of these
vectors are primitive vectors in $N_{\wPP}$. Let $(e_0, \ldots, e_n)$ be
the standard basis for $\ZZ^{n+1}$ and let $(f_0, \ldots, f_k)$ be the
standard basis for $\ZZ^{k+1}$.  Without loss of generality, we consider
the case in which the image of $e_0$ or $f_0$ is not primitive in
$N_{\wPP}$.  Let
\[
    r_1 =   \left(
            \begin{smallmatrix}
                a_0 \\ \vdots \\ a_n \\ 0 \\ \vdots \\ 0
            \end{smallmatrix}
            \right)
    \quad \text{and} \quad 
    r_2 =   \left(
            \begin{smallmatrix}
                -d_0 \\ \vdots \\ -d_n \\ 1 \\ \vdots \\ 1
            \end{smallmatrix}
            \right)
\]
be the columns of the matrix $T$.  Since the image of $e_0$ is not
primitive in $N_{\PP}$ or the image of $f_0$ is not, there exists a
vector $v \in \ZZ^{n+1} \oplus \ZZ^{k+1}$ such that 
\[
    e_0 + nv \in \Span\{r_1, r_2\} \quad \text{or} \quad f_0 + nv  \in
    \Span\{r_1, r_2\},
\]
for some integer $n$ such that $|n| > 1$.

In case the image of $e_0$ is not primitive in $N_{\wPP}$ then
\[
    e_0 + nv = c_1 r_1 + c_2 r_2 
\] 
for some integers $c_1$ and $c_2$.  Note that $n$ does not divide
$\gcd(c_1, c_2)$.  If we consider the components of $nv$, we find that 
$n \mid c_2$ and hence $n \mid c_1 a_i$ for all $i \neq 0$, and this
contradicts \eqref{Eq:wform}. Thus $n = 1$ and $v$ is primitive.  

In case the image of $f_0$ is not primitive in $N_{\wPP}$ then
\[
    f_0 + nv = c_1 r_1 + c_2 r_2 
\]
for some integers $c_1$ and $c_2$.  Considering the components of $nv$
we find that $n \mid c_2$ and $n \mid (c_2 - 1)$, hence $n = 1$, and
again we have a contradiction.
\end{proof}

We are now ready to write down the homogeneous coordinate ring and
irrelevant ideal $B$ for a toric $\PP(a_0, \ldots, a_n)$-bundle over
$\PP^k$.  

\begin{corollary}\label{C:coxring}
Let $(a_0, \ldots, a_n)$ be a sequence of positive integers satisfying
\eqref{Eq:wform}, and let ${\wPP}$ be a locally trivial toric
fibration over $\PP^k$ with fiber $\PP(a_0, \ldots, a_n)$. Then the
homogeneous coordinate ring of ${\wPP}$ is isomorphic to the
$\ZZ^2$-graded polynomial ring
\[
    \CC[X_0, \ldots, X_n, S_0, \ldots, S_k]
\]
where
\begin{equation*}
\begin{aligned}
	\deg X_i &= (a_i,-d_i) \quad \text{and}\\
	\deg S_j &= (0,1). 
\end{aligned}
\end{equation*}
Furthermore, the irrelevant ideal $B \subseteq S_X$ is given by 
\[
    B = (X_0 S_0, \ldots, X_n S_0, X_0 S_1, \ldots, X_n S_1, \ldots, X_0
    S_k, \ldots, X_n S_k).
\]
\end{corollary}

\begin{proof}
This follows immediately from Theorem~\ref{T:wpsbund/ps}, 
Proposition~\ref{P:wpstes}, and the construction in \cite{dC95} of the
homogeneous coordinate ring of a toric variety. 
\end{proof}

\section{Weighted projective space bundles as $\Proj$ of weighted
symmetric algebras}\label{S:proj}

In this section we view weighted projective space bundles over
projective space as $\Proj$ of weighted symmetric algebras of sums of
invertible sheaves and show that they are isomorphic to the toric
varieties in Construction~\ref{Co:toricwpb}.  For convenience, we adopt
the following definition.

\begin{definition}\label{D:wsheaf}
Let $X$ be a scheme, and let $(a_0, \ldots, a_n)$ be a
sequence of positive integers.  We define a \emph{weighted locally free
sheaf} with weights $(a_0, \ldots, a_n)$ to be a locally free sheaf of
$\Oh_X$-modules $\E$ together with an ordered decomposition $\E \cong
\E_0 \oplus \cdots \oplus \E_n$ such that $\E_i$ is an invertible sheaf
and such that the direct sum is to be interpreted as a graded sheaf with
$\E_i$ placed in degree $a_i$ for $0 \leq i \leq n$.  The weights will be
used in Definition~\ref{D:wpbund} when we form symmetric algebras of
weighted locally free sheaves. 
\end{definition}

\begin{remark}\label{R:wsheaf}
It certainly is possible not to require that the locally free sheaves
$\E_i$ in the decomposition in Definition~\ref{D:wsheaf} have rank one,
but we will restrict to this case as this is case that relates to toric
varieties.
\end{remark}

\begin{definition}\label{D:wpbund}
Let $X$ be a scheme.  Given a weighted locally free sheaf
$\E$ with weights $(a_0, \ldots, a_n)$ let $\mc{S}$ denote the weighted
symmetric algebra of $\E$ where we insist that $\E_i$ have homogeneous
degree $a_i$ in $\mc{S}$.  We define the \emph{weighted projective
bundle} associated to $\E$ to be the $X$-scheme
\[
\xymatrix{
    {\wPP(\E) := \Proj \mc{S}} \ar[r] & {X}.
}
\]
\end{definition}

We have the following two basic lemmas for weighted projective bundles.

\begin{lemma}\label{L:isloctriv}
Let $X$ be a variety over $\CC$, and let $\E$ be a weighted locally free
sheaf on $X$ with weights $(a_0, \ldots, a_n)$.  Then the weighted
projective bundle $\wPP(\E)$ is a locally trivial fiber bundle over $X$
with fiber the weighted projective space $\PP(a_0, \ldots, a_n)$.
\end{lemma}

\begin{proof}
This follows from the construction of $\Proj$ of a sheaf of graded
algebras as described, for example, in \cite[\S3]{egaII}.
\end{proof}

\begin{lemma}\label{L:twist}
Let $X$ be a scheme, let $\E$ be a weighted locally free
sheaf on $X$ with weights $(a_0, \ldots, a_n)$, and let $\mc{L}$ be an
invertible sheaf of $\Oh_X$-modules. Then there is a canonical
$X$-isomorphism
\[
    \wPP(\E) \cong \wPP \left( \bigoplus_{i=0}^{n} \E_i \otimes
    \mc{L}^{\otimes a_i} \right).
\]
\end{lemma}

\begin{proof}
This is just the application to our special case of
\cite[Proposition~3.1.8 (iii)]{egaII}.
\end{proof}

We also state the relative version of \cite[Lemma 5.7]{aI00}.

\begin{lemma}\label{L:wfprojsheaf}
Let $(a_0, \ldots, a_n)$ be a sequence of positive integers having no
common factor, let $X$ be a scheme, and let $\E =
\bigoplus_{i=0}^n \E_i$ be a weighted locally free sheaf of
$\Oh_X$-modules with weights $(a_0, \ldots, a_n)$.  Let $q$ be a
positive integer such that $q \mid a_i$ for $i > 0$, and define a new
weighted locally free sheaf $\E'$ with weights $(a_0, \frac{a_1}{q},
\ldots, \frac{a_n}{q})$ via
\[
    \E' := \E_0^{\otimes q} \oplus \bigoplus_{i = 1}^n \E_i.
\]
Then $\wPP(\E) \cong \wPP(\E')$.
\end{lemma}

We conclude this section by relating the two notions of weighted
projective bundles that we have described.  We point out that
Proposition~\ref{P:toricisproj} is a generalization of the exercise
\cite[Exercise p.  42]{wF93}.   

\begin{proposition}\label{P:toricisproj}
Let $(a_0, \ldots, a_n)$ be a sequence of positive integers satisfying
\eqref{Eq:wform}, and let $d_0, \ldots, d_n$ be integers.  Let $\E$
denote the weighted locally free sheaf 
\[      
    \E := \Oh(d_0) \oplus \cdots \oplus \Oh(d_n)
\]
on $\PP^k$ with weights $(a_0, \ldots, a_n)$.  Then $\wPP(\E)$ is
isomorphic to $\wPP$, the toric weighted projective bundle from
Construction~\ref{Co:toricwpb}.
\end{proposition}

\begin{proof}
We first fix some notation.  Write $\PP^k = \Proj \CC[S_0, \ldots S_k]$
and let $U_i = \Spec \CC[S_0/S_i, \ldots, S_k/S_i]$.  Let
$\{X_i\}_{i=0}^n$ be global coordinates on $\wPP(\E)$; this means that 
\[
    \wPP(\E)|_{U_i} \cong \Proj \left( \CC\left[\dfrac{S_0}{S_i}, \cdots,
    \dfrac{S_k}{S_i} \right][S_i^{d_0} X_0, \ldots, S_i^{d_n}
    X_n]\right)
\]
where $\deg S_i^{d_j} X_j = a_j$.  We find for $0 \leq i \leq k$ that
$\wPP(\E)|_{U_i}$ is isomorphic to the toric variety $\mb{A}^k \times
\PP(a_0, \ldots, a_n)$.  Let $U := \bigcap_i U_i$. Distinguishing $S_0$,
we find
\[
    \wPP(\E)|_U \cong \Proj \left( \CC \left[ 
    \left( \dfrac{S_1}{S_0}\right)^{\pm 1}, \cdots,
    \left(\dfrac{S_k}{S_0}\right)^{\pm 1} \right] [S_0^{d_0} X_0,
    \ldots, S_0^{d_n} X_n] \right)
\]
where $\deg S_0^{d_j} X_j = a_j$.  Distinguishing $S_0$
amounts to choosing a basis for the torus $\mb{T}^k$ in $\PP^k$ and
thereby identifying the lattice $M_{\PP^k}$ with $\ZZ^k$.  Similarly,
the coordinates $S_0^{d_i} X_i$ correspond to a basis for $\ZZ^{n+1}$ in
the context of the exact sequence \eqref{Eq:wpses}.  Therefore, the
inclusion
\[
\xymatrix{
    {\wPP(\E)|_U} \ar[r] & {\wPP(\E)|_{U_i}}}
\]
induces an isomorphism 
\[
\xymatrix{
    {M_{U_i} \oplus \ZZ^{n+1}} \ar[r]^-*{\phi_i} & {M_{\TT^k} \oplus
    \ZZ^{n+1}}},
\]
which is the identity if $i = 0$ and if $i > 0$ is given by the matrix
\begin{equation}\label{Eq:latmat}
\left(
\setcounter{MaxMatrixCols}{12}
\begin{matrix}
    0  &  1 & \hdots &  0 &  0 & \hdots & 0  & 0   & 0   & \hdots & 
        0 \\
    \vdots & \vdots & \ddots & \vdots & \vdots & \ddots &
        \vdots & \vdots & \vdots & \ddots & \vdots \\
    0  & 0  & \hdots &  1 &  0 & \hdots & 0  & 0   & 0   & \hdots & 
        0 \\
    -1 & -1 & \hdots & -1 & -1 & \hdots & -1 & d_0 & d_1 & \hdots &
        d_n \\
    0  & 0  & \hdots &  0 &  1 & \hdots & 0  & 0   & 0   & \hdots &
        0 \\
    \vdots & \vdots & \ddots & \vdots & \vdots & \ddots &
        \vdots & \vdots & \vdots & \ddots & \vdots \\
    0  & 0  & \hdots &  0 &  0 & \hdots & 1  & 0   & 0   & \hdots & 
        0 \\
    0  & 0  & \hdots &  0 &  0 & \hdots & 0  & 1   & 0   & \hdots & 
        0 \\
    0  & 0  & \hdots &  0 &  0 & \hdots & 0  & 0   & 1   & \hdots & 
        0 \\
    \vdots & \vdots & \ddots & \vdots & \vdots & \ddots &
        \vdots & \vdots & \vdots & \ddots & \vdots \\
    0  & 0  & \hdots &  0 &  0 & \hdots & 0  & 0   & 0   & \hdots & 
        1 
\end{matrix}
\right)
\end{equation}
where
\[
\left(
\begin{matrix}
    -1 & -1 & \hdots & -1 & -1 & \hdots & -1 & d_0 & d_1 & \hdots &
        d_n \\
\end{matrix}
\right)
\]
is the $i$-th row.  Taking the transpose of $\phi_i$ and composing with
the map in \eqref{Eq:wpses} we obtain isomorphisms 
\[
\xymatrix{
    {N_{\TT^k} \oplus N_{\PP(a_0, \ldots, a_n)}} \ar[r]^-*{\psi_i} &
    {N_{U_i} \oplus N_{\PP(a_0, \ldots, a_n)}}}.
\]
We put 
\[
    N_{\wPP(\E)} := N_{\TT^k} \oplus N_{\PP(a_0, \dots, a_n)},
\]
and the rays of the fan are the inverse images under the isomorphisms
$\psi_i$ of the rays of the fans for $U_i \times \PP(a_0, \ldots, a_n)$.
Observe that the inverse of the matrix \eqref{Eq:latmat} is given
by
\begin{equation}\label{Eq:latmatinv}
\left(
\setcounter{MaxMatrixCols}{12}
\begin{matrix}
    -1  &  -1 & \hdots &  -1 &  -1 &  -1 & \hdots & -1  & d_0   & d_1
    & \hdots & d_n \\
    1   &  0  & \hdots & 0   & 0   &   0 & \hdots & 0   & 0     & 0
    & \hdots & 0 \\
    0   &  1  & \hdots & 0   & 0   &   0 & \hdots & 0   & 0     & 0
    & \hdots & 0 \\
    \vdots & \vdots & \ddots & \vdots & \vdots & \vdots & \ddots &
    \vdots & \vdots & \vdots & \ddots & \vdots \\
    0   &  0  & \hdots & 1   & 0   &   0 & \hdots & 0   & 0     & 0
    & \hdots & 0 \\
    0   &  0  & \hdots & 0   & 0   &   1 & \hdots & 0   & 0     & 0
    & \hdots & 0 \\
    \vdots & \vdots & \ddots & \vdots & \vdots & \vdots & \ddots &
    \vdots & \vdots & \vdots & \ddots & \vdots \\
    0   &  0  & \hdots & 0   & 0   &   0 & \hdots & 1   & 0     & 0
    & \hdots & 0 \\
    0   &  0  & \hdots & 0   & 0   &   0 & \hdots & 0   & 1     & 0
    & \hdots & 0 \\
    0   &  0  & \hdots & 0   & 0   &   0 & \hdots & 0   & 0     & 1
    & \hdots & 0 \\
    \vdots & \vdots & \ddots & \vdots & \vdots & \vdots & \ddots &
    \vdots & \vdots & \vdots & \ddots & \vdots \\
    0   &  0  & \hdots & 0   & 0   &   0 & \hdots & 0   & 0     & 0
    & \hdots & 1 
\end{matrix}
\right)
\end{equation}
where the $i$-th column is the column whose only nonzero entry is $-1$.
Therefore the set of rays is the set of images in $N_{\wPP(\E)}$ of
the set of row vectors of the matrices \eqref{Eq:latmatinv} for all $i$.
They are given by the row vectors (as opposed to the column vectors)
because the matrix for the map $\psi_i$ is the transpose of the matrix
\eqref{Eq:latmatinv}.  It is now straightforward to check the relations
among this set of row vectors and verify that the fan for $\wPP(\E)$ is
isomorphic to the fan described in Construction~\ref{Co:toricwpb}.
\end{proof}

\section{Intersection numbers on weighted $\PP^1$-bundles over
$\PP^1$}\label{S:intth}

Let $a_0$ and $a_1$ be integers that are relatively prime.  In this
section we consider divisors on $\PP(a_0, a_1)$-bundles over $\PP^1$.
By Lemma~\ref{L:wfprojsheaf}, such bundles are isomorphic to
non-weighted $\PP^1$-bundles over $\PP^1$, and we check how divisors
transform under this isomorphism so that we may apply the formula
\eqref{Eq:intnumb} to obtain intersection numbers on
$\PP(a_0, a_1)$.  In order to fix terminology and notation, we recall
below the intersection theory of $\PP^1$-bundles over $\PP^1$.  See
\cite{wF98} for details.

\begin{example}\label{E:chowp3/p1}
Let $a$ and $b$ be integers, and form the locally free sheaf
\[
    \E := \Oh_{\PP^1}(a) \oplus \Oh_{\PP^1}(b)
\]
over $\PP^1$.  Let $H \in A^1\PP(\E)$ be the class of a fiber over
$\PP^1$, and let $L \in A^1\PP(\E)$ be the divisor class corresponding
to $\Oh_{\PP(\E)}(1)$.  Then the Chow ring of $\PP(\E)$ is given by
\begin{equation}\label{Eq:chowp3/p1}
    A^* \PP(\E) \cong \frac{\ZZ[H,L]}{(H^2, L^2 - (a+b)HL)}.
\end{equation}
Let $D_1$ be a divisor of type $(d_1,e_1)$, that is $D_1 \sim d_1L +
e_1H$, and let $D_2$ be a divisor of type $(d_2, e_2)$.  Then we find
from \eqref{Eq:chowp3/p1}, that the intersection number $D_1\cdot D_2$ is
given by
\begin{equation}\label{Eq:intnumb}
    D_1\cdot D_2 = \deg (D_1D_2) = (a+b)d_1d_2 + d_1e_2 + d_2e_1.
\end{equation}
The calculation in \eqref{Eq:intnumb} uses the fact that the cycle class
$HL$ is the class of a point in $A^2\PP(\E)$.
\end{example}

Let $d_0$ and $d_1$ be integers, and let $ \E := \Oh_{\PP^1}(d_0) \oplus
\Oh_{\PP^1}(d_1)$ be a weighted locally free sheaf with relatively prime
weights $(a_0, a_1)$.  Let $X_0$ and $X_1$ be global coordinates on
$\wPP(\E)$ as in the proof of Proposition~\ref{P:toricisproj}, and let
$S_0$ and $S_1$ be coordinates on the base $\PP^1$.  Each coordinate
$X_i$ defines an effective divisor $Z(X_i)$ on $\wPP(\E)$ by taking
$\Proj$ of the surjective map of weighted symmetric algebras induced by
the surjective map of sheaves 
\[
\xymatrix{
    {\E} \ar[r] & {\E/\Oh_{\PP^1}(d_i)}
}.
\]
The coordinates $S_i$ give divisors that are the pullbacks of
divisors from the base $\PP^1$.  From Lemma~\ref{L:wfprojsheaf} we know
that 
\[
    \wPP(\E) \cong \PP(\E'),
\] 
where $\E' = \Oh_{\PP^1}(a_1d_0) \oplus \Oh_{\PP^1}(a_0d_1)$.  The next
lemma shows how to realize a divisor given by global coordinates on
$\wPP(\E)$ as a divisor given by global coordinates on $\PP(\E)$. For
any real number $r$, the symbol $\lceil r \rceil$ denotes the smallest
integer greater than or equal to $r$.   

\begin{lemma}\label{L:globaldiv}
Let $\E = \Oh_{\PP^1}(d_0) \oplus \Oh_{\PP^1}(d_1)$ be a weighted
locally free sheaf with relatively prime weights $(a_0, a_1)$.  Let $D =
Z(S_0^{e_0}S_1^{e_1}X_0^{f_0}X_1^{f_1})$ be an effective divisor on
$\wPP(\E)$ defined by global coordinates.  Let $D'$ be the corresponding
divisor on the isomorphic variety $\PP(\E')$ where
\[
    \E' = \Oh_{\PP^1}(a_1d_0) \oplus \Oh_{\PP^1}(a_0d_1).
\]
Then $D' =
Z(S_0^{e_0}S_1^{e_1}X_0^{f'_0}X_1^{f'_1})$ where for $i = 0$ or $i=1$
and $j \neq i$ the integer $f'_i$ is given by
\begin{equation*}\label{Eq:fprime}
    f'_i =
        \left\lceil \frac{f_i}{a_j} \right\rceil. 
\end{equation*}
Furthermore, the divisor
$D'$ has type 
\[
    \left(
    \left\lceil\frac{f_0}{a_1}\right\rceil +
    \left\lceil\frac{f_1}{a_0}\right\rceil,
    e_0 + e_1 -
    \left\lceil\frac{f_0}{a_1}\right\rceil a_1d_0 - 
    \left\lceil\frac{f_1}{a_0}\right\rceil a_0d_1
    \right).
\]
\end{lemma}

\begin{proof}
The isomorphism of Lemma~\ref{L:wfprojsheaf} can be realized by taking
$\Proj$ of the degree $a_0a_1$ inclusion of sheaves of graded algebras
\[
\xymatrix{
    {\Sym_{\Oh_{\PP^1}}\E'} \ar@{^(->}[r]^-*{h}
    & {\widetilde{\Sym}_{\Oh_{\PP^1}}\E},
}
\]
where we have decorated the second ``$\Sym$'' to emphasize that it is a
weighted symmetric algebra.  The divisor $D$ is given by a subsheaf of
graded ideals
\[
    \mc{I}_0 \cdot \mc{J} \subseteq {\widetilde{\Sym}_{\Oh_{\PP^1}}\E},
\]
where $\mc{I}_0$ comes from the $S_i$ and is generated in degree $0$,
and where the degree $d$ part of $\mc{J}$ is given by 
\[
    \mc{J}_d = \left(
    \bigoplus_{\substack{k,l \geq 0,\\
                            a_0(f_0+k) + a_1(f_1+l) = d}}
    \Oh_{\PP^1}((f_0+k)d_0 + (f_1+l)d_1)\cdot X_0^{f_0+k}X_1^{f_1+l}
    \right).
\]
Now, the divisor $D'$ is defined by the sheaf of graded ideals 
\[
    h^{-1}(\mc{I}_0\cdot \mc{J}) \subseteq \Sym_{\Oh_{\PP^1}}\E', 
\]
and the first assertion follows easily from the description of
$\mc{J}_d$ above and the fact that $a_0$ and $a_1$ are relatively prime.

The second assertion follows from the fact that on any projective
bundle 
\[
    \PP(\Oh_{\PP^1}(d_0) \oplus \Oh_{\PP^1}(d_1)),
\]
the global coordinate $X_i$ defines a divisor of type $(1, -d_i)$, and
the global coordinate $S_j$ defines a divisor of type $(0,1)$.  This may
be checked directly, and it also follows from our description of the
bundle as a toric variety.
\end{proof}

\section{Linear systems}\label{S:linsys}

Before we address the question of finding anti-canonical hypersurfaces
in weighted projective space bundles, we give some results regarding
linear systems that we will need in what follows.  These results
are generalizations of the standard results for linear systems on
nonsingular varieties as discussed, for example, in \cite{rH77} or
\cite{pG94}.

Let $X$ be a normal, projective variety, and let $\mc{L}$ be a coherent
sheaf of $\Oh_X$-modules.  Recall that $\mc{L}$ is \emph{reflexive} if
the canonical map 
\[
\xymatrix{
    {\mc{L}} \ar[r] & {\mc{L}^{**}}
}
\]
is an isomorphism, where we define
\[
    \mc{L}^* := \sheafHom_{\Oh_X} (\mc{L}, \Oh_X).
\]
Furthermore, if $\mc{L}$ is a reflexive subsheaf of the function field
$k(X)$, then we say that $\mc{L}$ is \emph{divisorial}.  See \cite{mR80}
for details about divisorial sheaves.  The main point is that to any
Weil divisor $D$ on $X$, we may associate a divisorial sheaf $\Oh_X(D)$
and this assignment induces a one-to-one correspondence between the set
of Weil divisors on $X$ modulo linear equivalence and the set of
isomorphism classes of divisorial sheaves of $\Oh_X$-modules.  The
following fact is straightforward to check.

\begin{lemma}\label{L:linsys}
Let $D$ be a Weil divisor on $X$.  The set of divisors that are linearly
equivalent to $D$ is in one-to-one correspondence with the set of closed
points of the projective space
\[
    \PP H^0(X, \Oh_X(D)).
\]
\end{lemma}

\begin{definition}\label{D:linsys}
Let $X$ be a normal projective variety, and let $D$ be a Weil divisor
on $X$.  The \emph{complete linear system}, denoted $|D|$, is the set of
all Weil divisors that are linearly equivalent to $D$.  A \emph{linear
system} is a linear subspace of a complete linear system, where we
regard the complete linear system as having the structure of a
projective space guaranteed by Fact~\ref{L:linsys}.
\end{definition}

\begin{definition}\label{D:bs}
The \emph{base locus} of a linear system $L$ is the set-theoretic
intersection of the supports of the members of $L$.
\end{definition}

\begin{definition}\label{D:V_L}
Let $L \subseteq |D|$ be a linear system for some Weil divisor $D$ on
$X$.  We denote by $V_L$ the associated linear subspace of
$H^0(X,\Oh_X(D))$.
\end{definition}

Now suppose that $X_{\Sigma}$ is a simplicial toric variety associated
to some fan $\Sigma \subseteq N_{\RR}$ for some lattice $N$.  Let $S$ be
the homogeneous coordinate ring of $X_{\Sigma}$.  Recall that $S$ is an
$A^1 X_{\Sigma}$-graded polynomial ring.  Let $D$ be a Weil divisor on
$X$.  From \cite[Proposition~1.1]{dC95} we find that   
\begin{equation}\label{Eq:sections}
    H^0(X_\Sigma, \Oh_{X_\Sigma}(D)) \cong S_{[D]},
\end{equation}
where $[D]$ denotes the class of the divisor $D$ in $A^1(X_\Sigma)$, and
$S_{[D]}$ denotes the degree $[D]$ part of the graded ring $S$.  One
easily checks that the zero scheme (see \cite{dC95}) of a homogeneous
form $F \in S_{[D]}$ is a divisor on $X$ in the class of $D$ and, hence,
that linear systems may be described as linear subspaces of $S_{[D]}$.
We illustrate the case for weighted projective bundles in the following
example.

\begin{example}\label{R:genlinsys}
Any complete linear system $|D|$ on $\wPP(\E)$ may be written as the zero
locus of forms
\begin{equation}\label{Eq:linsysforms}
    \sum \phi_{e_0 \cdots e_n}(S_0,\ldots,S_k)X_0^{e_0}\dots X_n^{e_n}
    \subseteq S_{[D]}
\end{equation}
where the exponents $e_i$ satisfy $\sum a_ie_i = C$ for
some integer $C$, and $\phi_{e_0 \cdots e_n}$ is any form in
$H^0(\PP^k,\Oh(\eta))$ for some integer-valued function $\eta$ that
depends on the exponents $e_i$.  In case $|D| = -K$, the function $\eta$
is given by 
\begin{equation}\label{Eq:eta}
    \eta(e_0, \ldots, e_n) := \sum d_i(e_i - 1) + k + 1.
\end{equation}
\end{example}

We now check that the property of a linear system on a toric variety
having a quasi-smooth member is an open condition.  

\begin{lemma}\label{L:opencond}
Let $X$ be a simplicial toric variety that is proper over $\CC$, and let
$L$ be a nonempty linear system on $X$.  Then there exists a
quasi-smooth member of $L$ if and only if the general member of $L$ is
quasi-smooth.
\end{lemma}

\begin{proof}
One direction is clear.  Since $L$ is nonempty and the general member is
quasi-smooth, it must be the case that $L$ contains a quasi-smooth
member.  

We now establish the converse.  Assume that $L$ contains a quasi-smooth
member.  Consider the closed incidence subvariety $Z_{sing} \subseteq L
\times X$ given by 
\[
    Z_{sing} := \{[D] \times p \mid p \in \sing D\},
\]
where $\sing D$ denotes the locus where $D$ fails to be quasi-smooth.
Let $\pi_L$ be the projection of $Z_{sing}$ onto $L$.  The result now
follows from the properness of $\pi_L$.
\end{proof}
 
One of the main tools in the analysis of linear systems is Bertini's
Theorem.  We check that a version of Bertini's Theorem holds that
addresses the quasi-smoothness of members of linear systems on toric
varieties. 

\begin{proposition}\label{P:qbertini}
Let $X$ be a simplicial toric variety that is proper over $\CC$, and let
$L$ be a linear system on $X$.  Then the general member of $L$ is
quasi-smooth away from the base locus of $L$.
\end{proposition}

\begin{proof}
Let $S$ be the homogeneous coordinate ring of $X$, and let 
\[
\xymatrix{
    {\Spec S - Z(B)} \ar[r]^-*{q} & {X}
}
\]
be the quotient map (see Section~\ref{S:conventions}).  Since
quasi-smoothness of a divisor on $X$ is determined by the nonsingularity
of a hypersurface on $\Spec S - Z(B)$, the argument given in
the proof of Bertini's Theorem in \cite[p.~137]{pG94} applies in our
case to establish the proposition.
\end{proof}

It will be convenient to have the notion of the restriction of a linear
system to the fiber of a weighted projective space bundle.  

\begin{definition}\label{D:pblinsys}
Let $X$ be a simplicial toric variety, let $j:Y \hookrightarrow X$ be a
well-formed quasi-smooth subvariety, and let $L \subseteq |D|$ be a
linear system on $X$ for some Weil divisor $D$.  Since $Y$ is
well-formed, the sheaf $j^*\Oh_X(D)$ is divisorial.  We define the
\emph{restricted linear system}, denoted $L|_Y$ to be the linear system
associated to the image of $V_L$ under the natural $\CC$-linear map
\[
\xymatrix{
    {H^0(X,\Oh_X(D))} \ar[r] & {H^0(Y, j^*\Oh_X(D))}.
    }
\]
\end{definition}

\begin{example}\label{E:fibindep}
Let $|D|$ be a complete linear system on a weighted projective space
bundle, and let $F \cong \PP(a_0, \ldots, a_n)$ be a fiber of the
bundle.  It follows easily from the description given in
Example~\ref{R:genlinsys} that the restricted linear system $|D||_F$ is
independent of the choice of fiber $F$ when regarded as a linear system
on $\PP(a_0, \ldots, a_n)$.
\end{example}

\begin{proposition}\label{P:gensmoothfib}
Let $|D|$ be a complete linear system on $\wPP(\E)$ over $\PP^k$ for
some divisor $|D|$, and let $F \cong \PP(a_0, \ldots, a_n)$ be a fiber.
If the general member of $|D|$ is quasi-smooth, then so is the general
member of the restricted linear system ${L : =|D||_F}$.
\end{proposition}

\begin{proof}
Suppose, to the contrary, that every member of $L$ fails to be
quasi-smooth.  We will show that the same is true for $|D|$.  Let $G$ be
the equation for a member of $|D|$.  Restricting $G$ to any fiber,
the assumption that every member of $L$ fails to be quasi-smooth implies
that the system of equations
\[
    G = \frac{\del G}{\del X_0} = \cdots = \frac{\del G}{\del X_n} = 0
\]
has a nontrivial solution $P$.  Furthermore, by
Proposition~\ref{P:qbertini}, we may assume that $(X_0(P),\ldots,
X_n(P)) \in \Bs L$.  We claim that $P$ is also a solution to the system
of equations
\[
    \frac{\del G}{\del S_0} = \cdots = \frac{\del G}{\del S_k} = 0.
\]
Indeed, restricting the equations above to the fiber $F_P$ that
contains $P$ yields a system of $(k+1)$ equations, each of which defines
a member of $L$. But $P|_{F_P} \in \Bs L$ by assumption, and the claim
follows.  
\end{proof}

To conclude this section we state a result from \cite{aI00} that
addresses the problem of determining whether a linear system on
$\PP(a_0, \ldots, a_3)$ has a quasi-smooth member.  In \cite{aI00} the
result is stated only for complete linear systems, but the result is a
corollary of \cite[8.1~Theorem]{aI00} whose proof applies more generally
to linear systems whose associated vector space of homogeneous forms is
generated by monomials.  

\begin{definition}\label{D:monls}
Let $L$ be a linear system on weighted projective space.  If the
associated vector space of forms is generated by monomials, then we say
that $L$ is a \emph{monomial linear system}.
\end{definition}

\begin{proposition}[{\cite[8.5~Corollary]{aI00}}]\label{P:wp3qsmooth}
 Let $(a_0, \ldots, a_3)$ be a sequence of positive integers satisfying
\eqref{Eq:wform}, and let $L$ be a monomial linear system on
$\PP(a_0, \ldots, a_3)$ of degree $d$ where $d > a_i$ for all $0 \leq i
\leq 3$.  Then the general member of $L$ is quasi-smooth if and only if
the following two conditions hold. 
\begin{enumerate}
    \item For all $0 \leq i \leq 3$ there exists a monomial in $V_L$
    that does not involve $X_i$.
    \item For all $0 \leq i \leq 3$ there exists a monomial in $V_L$ of
    the form $X_i^pX_{e_i}$ for some $0 \leq e_i \leq 3$.
    \item For all $0 \leq i < j \leq 3$ 
    \begin{enumerate}
        \item there exists a monomial in $V_L$ of the form $X_i^pX_j^q$
        or
        \item there exist monomials in $V_L$ of the form
        $X_i^{p_1}X_j^{q_1}X_{e_1}$ and $X_i^{p_2}X_j^{q_2}X_{e_2}$ for
        distinct $e_i$.  
    \end{enumerate}
\end{enumerate}
\end{proposition}

\section{Well-formedness and anti-canonical
hypersurfaces}\label{S:wform}

This section is concerned with the condition of well-formedness as
given in Definition~\ref{D:wformsub}.  In Section~\ref{S:achs},
we will apply the adjunction formula to a quasi-smooth anti-canonical
hypersurface in a weighted projective bundle in order to show that it is
a Calabi-Yau variety.  For the adjunction formula to work as
expected, we will need to know that the hypersurface in question is
well-formed.  As we now establish, it turns out that well-formedness is
automatic for quasi-smooth anti-canonical hypersurfaces.  In \cite{aI00}
Iano-Fletcher considered the problem of well-formedness of hypersurfaces
in weighted projective space, and we take up this case first.  

In order to obtain information about the well-formedness of
hypersurfaces in weighted projective space, we need to understand the
singular locus of weighted projective space.  All the singularities are
finite cyclic quotient singularities, and we recall notation
for such singularities below.

\begin{definition}\label{D:quotsing}
Let $a_1, \ldots, a_n$ be integers, let $r$ be a positive integer, and
let $x_1, \ldots, x_n$ be coordinates on $\mb{A}^n$.  Let $\ZZ_r$ act on
$\mb{A}^n$ via 
\begin{equation}\label{Eq:action}
    x_i \mapsto \varepsilon^{a_i} x_i \quad \text{for all $i$},
\end{equation}
where $\varepsilon$ is a fixed primitive $r$-th root of unity.  Let $X$
be an algebraic variety over $\CC$.  A singularity $Q \in X$ is a
\emph{(quotient) singularity of type $\frac{1}{r}(a_1, \ldots, a_n)$} if
$(X,Q)$ is locally analytically isomorphic to the quotient
$(\mb{A}^n/\ZZ_r,O)$ under the action defined in \eqref{Eq:action},
where $O$ denotes the image of the origin under the quotient map.
\end{definition}

\begin{lemma}[{\cite[\S5.15]{aI00}}]\label{L:singwps}
Let $(a_0, \ldots, a_n)$ be a sequence of integers satisfying
\eqref{Eq:wform}.  Let $Z_{i_1\cdots i_d} \subseteq \PP(a_0,
\ldots, a_n)$ be the closed subvariety defined by 
\[
    X_{i_1} = \cdots = X_{i_d} = 0
\]
for distinct $i_1, \ldots, i_d$.  The general point of $Z_{i_1\cdots
i_d}$ is locally analytically isomorphic to $(O,Q) \in \mb{A}^{n-d}
\times Y$ where $Q \in Y$ is a singularity of type 
\[
    \frac{1}{h_{i_1\cdots i_d}} (a_{i_1}\cdots, a_{i_d}),
\]
and where $h_{i_1\cdots i_d}$ is the greatest common divisor of the
complement of $(a_{i_1}, \ldots, a_{i_d})$ in $(a_0, \ldots, a_d)$.
\end{lemma}

The next lemma is a generalization of \cite[6.10]{aI00}, in which it is
stated (without proof) only for complete linear systems.  

\begin{lemma}[{\cite[\S6.10]{aI00}}]\label{L:wformhypwps}
Let $(a_0, \ldots, a_n)$ be a sequence of positive integers satisfying
\eqref{Eq:wform}, and let $L$ be a monomial linear system of degree
$d$ on $\PP(a_0, \ldots, a_n)$. The general member of $L$ is well-formed
if and only if for all $0 \leq i < j \leq n$
\begin{enumerate}
    \item  there exists a monomial $M \in V_L$ such that $M \notin (X_i,
    X_j)$, or
    \item we have $h_{ij} = 1$ where $h_{ij}$ is as defined in the
    statement of Lemma~\ref{L:singwps}.
\end{enumerate}
\end{lemma}

\begin{proof}
First suppose that the general member $L$ is well-formed and that
condition 1 fails to hold for some $0 \leq i < j \leq n$.  We must
show that condition 2 holds.  The failure of condition 1 to hold means
precisely that the subvariety $Z_{ij}$ of $\PP(a_0, \ldots, a_n)$
defined by $X_i = X_j = 0$ is contained in the base locus of $L$.  But
if the general member of $L$ is well-formed it must be the case that the
general point of $Z_{ij}$ is nonsingular.  According to
Lemma~\ref{L:singwps} we must have $h_{ij} =1$ and this is condition 2.  

Now assume that condition 1 holds for all $0 \leq i < j \leq n$.  Then
the codimension of the base locus of $L$ is greater than two, and it
follows that the general member $L$ is well-formed.  Suppose that
condition one fails to hold for some $0 \leq i < j \leq n$, but
condition two does hold.  Then the subvariety $Z_{ij}$ defined in the
previous paragraph is in the base locus, but its general point is
nonsingular.  Again, it follows that the general member of $L$ is
well-formed.  
\end{proof}

\begin{proposition}\label{P:anticanwform}
Let $(a_0, \ldots, a_3)$ be a sequence of positive integers satisfying
\eqref{Eq:wform}.  Let $L \subseteq |-K_{\PP(a_0, \ldots, a_3)}|$ be a
monomial linear system.  If the general member of $L$ is quasi-smooth,
then the general member of $L$ is well-formed. 
\end{proposition}

\begin{proof}
Suppose that condition 1 of Lemma~\ref{L:wformhypwps} fails to hold for
some $0 \leq i < j \leq 3$.  This means that the base locus of $L$ has
an irreducible component given by $X_i = X_j = 0$.  Let $\{i', j'\}$
denote the complement of $\{i, j\}$ in the set $\{0, \ldots, 3\}$.  Then
according to Proposition~\ref{P:wp3qsmooth}, there exist monomials in
$V_L$ of the from $X_{i'}^{p_1}X_{j'}^{q_1}X_{e_1}$ and
$X_{i'}^{p_2}X_{j'}^{q_2}X_{e_2}$ for distinct $e_1$ and $e_2$ such that
$\{e_1, e_2\} = \{i,j\}$  In particular, this means that 
\[
    p_1a_{i'} + p_2a_{j'} + a_{e_1} = \sum_{l=0}^3 a_l. 
\]
Since we are assuming that \eqref{Eq:wform} holds, the equation
above implies that 
\[
    \gcd(a_{i'},a_{j'}) = 1
\]
as required to satisfy condition 2 of Lemma~\ref{L:wformhypwps}.
Therefore, the general member of $L$ is well-formed.
\end{proof}

\begin{corollary}\label{C:anticanwform}
Let $\wPP(\E)$ be a weighted $\PP(a_0, \ldots, a_3)$-bundle over
$\PP^1$.  If the general member of $|-K_{\wPP(\E)}|$ is quasi-smooth,
then it is also well-formed.
\end{corollary}

\begin{proof}
If the general member of $|-K_{\wPP(\E)}|$ is not well-formed, then the
same is true for the linear system $L$ on $\PP(a_0, \ldots, a_3)$
obtained by restricting $|{-K_{\wPP(\E)}}|$ to a fiber.  The result now
follows easily from Proposition~\ref{P:gensmoothfib} and
Proposition~\ref{P:anticanwform}.
\end{proof}

\section{Quasi-smooth and well-formed anti-canonical
hypersurfaces}\label{S:achs}

We come now to our main result, in which we classify all $\PP(a_0,
\ldots, a_3)$-bundles over $\PP^1$ whose anti-canonical linear systems
have a quasi-smooth member.  We first check that quasi-smooth
anti-canonical hypersurfaces in weighted projective bundles over $\PP^k$
are necessarily Calabi-Yau varieties fibered by smaller dimensional
Calabi-Yau varieties.  Due to the various notions in the literature, we
include here the definition of a Calabi-Yau variety that we use.

\begin{definition}[{\cite[Definition~1.4.1]{dC99}}]\label{D:cy}
Let $X$ be a normal variety that is proper over $\CC$ and has at worst
canonical singularities.  We say that $X$ is a \emph{Calabi-Yau variety}
if $X$ has trivial dualizing sheaf and 
\[
    H^1(X, \Oh_X) = \dots = H^{\dim X - 1}(X, \Oh_X) = 0.
\] 
\end{definition}

In order to conclude that quasi-smooth anti-canonical hypersurfaces in
weighted projective space bundles are fibered by smaller dimensional
Calabi-Yau varieties, we must first check that quasi-smooth
anti-canonical hypersurfaces of weighted projective space are
themselves Calabi-Yau varieties.

\begin{proposition}\label{P:weightedcy}
Let $(a_0, \ldots, a_n)$ be a sequence of integers satisfying
\eqref{Eq:wform}.  A quasi-smooth and well-formed member  
\[
    X \in |-K_{\PP(a_0, \ldots, a_n)}|
\]
is a Calabi-Yau variety. 
\end{proposition}

\begin{proof}
First we show that the dualizing sheaf $\omega_X$ is trivial.  For any
well-formed, quasi-smooth hypersurface $Z_d$ of degree $d$ in $\PP(a_0,
\ldots, a_n)$ we have 
\[
    \omega_{Z_d} \cong \Oh_{Z_d}\left(d - \sum a_i\right)
\]
(see \cite[Theorem~3.3.4]{iD82} or \cite[\S6.14]{aI00}).  Since $X$ is
anti-canonical, it has degree $\sum a_i$, and it follows that $\omega_X
\cong \Oh_X$.  Since $X$ is quasi-smooth, it is an orbifold, and a
Gorenstein orbifold necessarily has canonical singularities (see
\cite[Proposition 1.7]{mR80}).  Finally, the necessary vanishing of
cohomology follows from the known cohomology of hypersurfaces in
weighted projective space.  (See \cite[Section~3.4.3]{iD82} or
\cite[7.1~Lemma]{aI00}).
\end{proof}

\begin{theorem}\label{T:iscy} 
Let $\wPP(\E)$ be a weighted projective bundle 
over $\PP^k$.  A quasi-smooth and well-formed member of
$|-K_{\wPP(\E)}|$ is a Calabi-Yau variety whose general fiber over
$\PP^k$ is a weighted Calabi-Yau hypersurface.
\end{theorem}

\begin{proof} 
Let $X \in |-K_{\wPP(\E)}|$ be quasi-smooth and well-formed.  Note that
if $X$ is quasi-smooth, it is automatically irreducible (this follows
since $X$ is transverse to a fiber).  We first show that $K_X = 0$.  Let
$X_0 = \wPP(\E)_{sm} \cap X$ where $\wPP(\E)_{sm}$ denotes the
nonsingular locus of $\wPP(\E)$, and let $i: X_0 \longrightarrow X$ be
the inclusion map.  Since $X$ is well-formed, it follows that $X_0$ is
smooth and that the codimension of $X \setminus X_0$ in $X$ is at least
two.  The adjunction formula applied to $X_0$ and $\wPP(\E)_{sm}$ yields 
\begin{align*}
	\Oh_X(K_X) &\cong i_* \Oh_{X_0}(K_{X_0}) \\
		&\cong i_*[(\Oh_{\wPP(\E)_{sm}}(K_{\wPP(\E)_{sm}})
			\otimes \Oh_{\wPP(\E)_{sm}}(X_0) \otimes \Oh_{X_0}] \\
		&\cong i_* \Oh_{X_0} \\
		&\cong \Oh_X.
\end{align*}
The last isomorphism follows from the fact that both sheaves are
isomorphic to the divisorial sheaf associated to the zero divisor.  The
dualizing sheaf of $X$ is given by $\Oh_X(K_X) \cong \Oh_X$, so $X$ is
Gorenstein, and a Gorenstein orbifold necessarily has canonical
singularities (see \cite[Proposition 1.7]{mR80}).    

We next show the necessary vanishing of cohomology.  We use the
following exact sequence.  
\begin{equation}\label{Eq:iscyles}
\xymatrix{
	{0} \ar[r] & {\Oh_{\wPP(\E)}(K_{\wPP(\E)})} \ar[r] &
		{\Oh_{\wPP(\E)}} \ar[r] & {\Oh_X} \ar[r] & 0 }
\end{equation}
Taking the long exact sequence in cohomology, we find it suffices
to show that $H^i(\wPP(\E), \Oh_{\wPP(\E)}) = 0$ for $1 \leq i \leq
k+n-2$ and $H^j(\wPP(\E), \Oh_{\wPP(\E)}(K_{\wPP(\E)})) = 0$ for $2 \leq
j \leq k+n-1$.  By Serre duality we have 
\[
	H^j(\wPP(\E), \Oh_{\wPP(\E)}(K_{\wPP(\E)})) \cong 
	H^{k+n-j}(\wPP(\E), \Oh_{\wPP(\E)})^*,
\]
and for any toric variety $Z$, we have $H^i(Z, \Oh_Z) = 0$ for $i > 0$
(\cite[Corollary 7.4]{vD78}); this completes the proof that $X$ is a
Calabi-Yau variety.  

For the second assertion, note that generic smoothness over $\PP^k$ of
$q^{-1}(X)$, where $q$ is the quotient map in
Definition~\ref{D:quasismooth}, implies that the general fiber of $X
\longrightarrow \PP^k$ is quasi-smooth.  Therefore, the general fiber is
a well-formed anti-canonical quasi-smooth hypersurface in $\PP(a_0,
\ldots, a_n)$, whence a weighted Calabi-Yau hypersurface by
Proposition~\ref{P:weightedcy}.  
\end{proof}

It remains to find necessary and sufficient conditions for the
anti-canonical system of a weighted projective bundle to have a
quasi-smooth member.  We begin with a lemma that allows us to assume
that the weighted sheaf from which we form the weighted projective
bundle has a certain normalized form.

\begin{lemma}\label{L:normform}
Every weighted projective space bundle over $\PP^k$ is isomorphic to
$\wPP(\E)$, where $\E$ is a weighted locally free sheaf of rank $n+1$
with weights $(a_0, \ldots, a_n)$ satisfying \eqref{Eq:wform} and
such that 
\[
    \E \cong \bigoplus_{i=0}^n \Oh_{\PP^k}(d_i)
\]
where 
\begin{equation}\label{Eq:normform1}
    0 \leq d_0 \leq d_1 \leq \cdots \leq d_n,
\end{equation}
and there exists an $0 \leq i_0 \leq n$ such that 
\begin{equation}\label{Eq:normform2}
    d_{i_0} < a_{i_0}.
\end{equation}
\end{lemma}

\begin{proof}
Let $\E$ be a weighted locally free sheaf with weights $(a_0, \ldots,
a_n)$.  Since we are only considering weighted locally free sheaves that are
sums of invertible sheaves (see Definition~\ref{D:wsheaf}), we may
assume that 
\[
    \E \cong \bigoplus_{i=0}^n \Oh_{\PP^k}(e_i)
\]
for some integers $e_0, \ldots, e_n$ and that $\E$ is weighted by $(a_0,
\ldots, a_n)$. Given any integer $l$ we find from Lemma~\ref{L:twist}
that $\wPP(\E) \cong \wPP(\E')$ where 
\[
    \E' = \bigoplus_{i=0}^n \Oh_{\PP^k}(e_i + la_i).
\]
Let $\lambda$ be the smallest integer such that $e_i + \lambda a_i \geq
0$ for all $0 \leq i \leq n$, define 
\[
    d_i := e_i + \lambda a_i,
\]
and renumber so that $d_0 \leq d_1 \leq \cdots \leq d_n$.  If $d_i \geq
a_i$ for all $0 \leq i \leq n$, then we have $e_i + (\lambda-1)a_i \geq 0$,
but this contradicts the minimality of the integer $\lambda$.  Therefore $d_0,
\ldots, d_n$ give the required integers.
\end{proof}

We are now in position to prove our classification result for
anti-canonical systems on $\PP(a_0, \ldots, a_3)$-bundles over $\PP^1$.

\begin{theorem}\label{T:qsmoothk3fib}
Let $(d_0, \ldots, d_3)$ be a sequence of integers satisfying
\eqref{Eq:normform1} and \eqref{Eq:normform2}, and let $(a_0, \ldots,
a_3)$ be a sequence of positive integers satisfying
\eqref{Eq:wform}.  Let $\E$ be the weighted locally free sheaf
\[
    \E := \Oh_{\PP^1}(d_0) \oplus \cdots \oplus \Oh_{\PP^1}(d_3)
\]
with weights $(a_0, \ldots, a_3)$.  Let $L$ denote the linear system on
$\PP(a_0, \ldots, a_3)$ obtained by restricting $|-K_{\wPP(\E)}|$ to a
fiber, let $V_L$ be its associated vector space of homogeneous
forms, and let $\eta$ be defined via \eqref{Eq:eta}.  Then the linear
system $|-K_{\wPP(\E)}|$ has a well-formed quasi-smooth member if and
only if the following conditions hold.
\begin{enumerate}
    \item The general member of the linear system $L$ is quasi-smooth.
    \item For all $0 \leq i \leq 3$
    \begin{enumerate}
        \item there exists a monomial of the form $X^p_i$ in $V_L$, or
        \item there exists an integer $j \neq i$ and a monomial $M \in
        V_L$ of the form $X_jX_i^{p}$ such that $\eta(M) = 0$, or
        \item there exist distinct integers $j_1 \neq i$ and $j_2 \neq
        i$ and two monomials $M_1$ and $M_2$ of the form
        $X_{j_1}X_{i}^{p}$ and $X_{j_2}X_{i}^{p}$ respectively such that
        $\eta(M_1) > 0$ and $\eta(M_2) > 0$.
    \end{enumerate}
    \item For all $0 \leq i < j \leq 3$
    \begin{enumerate}
        \item there exists a monomial of the form $X_i^pX_j^q \in V_L$
        or 
        \item there exist monomials of the form
        \[
        N_1 = X_i^{p_1}X_j^{q_1}X_{e_1} \in V_L
        \]
        and 
        \[
        N_2 = X_i^{p_2}X_j^{q_2}X_{e_2} \in V_L
        \]
        in $V_L$ such that $e_1$ and $e_2$ are distinct and
        such that
        \begin{multline}\label{Eq:qsmoothintnum}
        \left(
        \left\lceil\frac{p_1}{a_j}\right\rceil +
        \left\lceil\frac{q_1}{a_i}\right\rceil
        \right)\eta(N_2) +
        \left(
        \left\lceil\frac{p_2}{a_j}\right\rceil +
        \left\lceil\frac{q_2}{a_i}\right\rceil
        \right)\eta(N_1) +\\
        \left(
        \left\lceil\frac{p_1}{a_j}\right\rceil 
        \left\lceil\frac{p_2}{a_j}\right\rceil +
        \left\lceil\frac{q_1}{a_i}\right\rceil
        \left\lceil\frac{q_2}{a_i}\right\rceil
        \right)(a_id_j - a_jd_i) = 0.
        \end{multline}
    \end{enumerate}
\end{enumerate}
\end{theorem}

\begin{proof} Recall that if the general member of $|-K_{\wPP(\E)}|$ is
quasi-smooth, then it is automatically well-formed by
Corollary~\ref{C:anticanwform}.  Now suppose that the general member of
$|-K_{\wPP(\E)}|$ is quasi-smooth.  From
Proposition~\ref{P:gensmoothfib} we see that condition 1 holds.  We must
show that conditions 2 and 3 hold.  We proceed according to the
codimension of the base locus of $|-K_{\wPP(\E)}|$.  The equations
that define the base locus are given by monomials in the variables $X_0,
\ldots, X_3$.  If the base locus is empty, then conditions 2.(a) and
3.(a) hold.  Also, there can be no fixed component because then $L$
would have a fixed component, and this would contradict the fact that
condition 1 holds.  This leaves two possibilities: either the base locus
of $|-K_{\wPP(\E)}|$ has codimension three or it has codimension two.  

We first consider the case in which 
\[
    \codim \Bs |-K_{\wPP(\E)}| = 3,
\]
and we suppose that $C$ is a component of the base locus given by 
\[
    X_{j_1} = X_{j_2} = X_{j_3} = 0
\]
for 
\[
    0 \leq j_1 < j_2 < j_3 \leq 3.
\]
Let $F$ denote the equation for a general member of $|-K_{\wPP(\E)}|$.
Since $Z(F)$ is quasi-smooth we find, in particular, that the system of
equations
\begin{equation}\label{Eq:k3fibbs0}
    \frac{\del F}{\del X_{j_1}} = 
    \frac{\del F}{\del X_{j_2}} = 
    \frac{\del F}{\del X_{j_3}} = 0
\end{equation}
has no nontrivial solution on $C$.  Restricting to $C$, the equations in
\eqref{Eq:k3fibbs0} become
\begin{equation}\label{Eq:k3fibbs0res}
    \phi_{j_1}X_i^{p_1} =
    \phi_{j_2}X_i^{p_2} =
    \phi_{j_3}X_i^{p_3} = 0,
\end{equation}
for $i \notin \{j_1,j_2,j_3\}$ and for forms 
\[
    \phi_{j_l} \in 
    H^0(\PP^1,\Oh_{\PP^1}(\eta(X_{j_l}X_i^{p_l}))).
\]
At least one of the forms must be nonzero otherwise there could be no
quasi-smooth member.  If only one form is nonzero, it must have degree
zero, otherwise it would have a zero.  This is condition 2.(b).  If at
least two of the $\phi_{jl}$ are nonzero, and neither has degree zero,
then condition 2.(c) holds.  Finally, since the dimension of the base
locus is zero, it follows that condition 3.(a) holds. 

Now assume that 
\[
    \codim \Bs |-K_{\wPP(\E)}| = 2.
\]    
If the base locus has a component of dimension zero, then we are in the
situation of the previous paragraph.  Let $C$ be a component of the base
locus given by $X_{e_1} = X_{e_2} = 0$ and let $\{i,j\}$ denote the
complement of $\{e_1,e_2\}$ in $\{0,1,2,3\}$.  Observe that 
\[
    C \cong \wPP(\Oh_{\PP^1}(d_i) \oplus \Oh_{\PP^1}(d_j))
\]
with weights $(a_i, a_j)$.  Furthermore, we claim that $\gcd(a_i, a_j) =
1$.  Indeed, since condition one is satisfied,
Proposition~\ref{P:wp3qsmooth} ensures the existence of a monomial of
the form $X_i^{p_1}X_j^{q_1}X_{e_1}$. Hence we have 
\[
    p_1a_i + p_2a_j + a_{e_1} = \sum a_i,
\]
from which the claim follows since the weights $(a_0, \ldots, a_3)$
satisfy \eqref{Eq:wform}.  Therefore, it follows from
Lemma~\ref{L:globaldiv} that 
\[
    C \cong \PP(\Oh_{\PP^1}(a_jd_i) \oplus \Oh_{\PP^1}(a_id_j)).
\]

Arguing as before we find that the quasi-smoothness of the general
member $Z(F)$ implies that the system of equations
\begin{equation}\label{Eq:k3fibbs1}
    \frac{\del F}{\del X_{e_1}} = \frac{\del F}{\del X_{e_2}} = 0
\end{equation}
has no nontrivial solution on $C$.  If we restrict the partial
derivatives \eqref{Eq:k3fibbs1} to $C$ we obtain the equation for two
divisors 
\[
    D_{e_1} := Z \left(\left.\frac{\del F}{\del X_{e_1}}\right|_C\right) 
    \quad \text{and} \quad 
    D_{e_2} := Z \left(\left.\frac{\del F}{\del X_{e_2}}\right|_C\right)
\]
on $C$ that meet transversally.  Therefore, the existence of a
nontrivial solution on $C$ to \eqref{Eq:k3fibbs1} is equivalent to the
intersection number $D_{e_1}\cdot D_{e_2}$ not vanishing.  Furthermore,
since condition 1 is satisfied, Proposition~\ref{P:wp3qsmooth}
guarantees the existence of monomials of the form $N_1$ and $N_2$ in
$V_L$, and we must check that \eqref{Eq:qsmoothintnum} is satisfied.
The presence of the monomials $N_1$ and $N_2$ in $V_L$ allows us to use
Lemma~\ref{L:globaldiv} to determine the type of the two divisors on
$C$.  Referring to Example~\ref{E:chowp3/p1} for notation and
terminology, we find that 
\[
    D_{e_1} \text{ has type } 
    \left(
    \left\lceil\frac{p_1}{a_j}\right\rceil + 
    \left\lceil\frac{q_1}{a_i}\right\rceil,
    \eta(N_1) -
    \left\lceil\frac{p_1}{a_j}\right\rceil a_jd_i -  
    \left\lceil\frac{q_1}{a_i}\right\rceil a_id_j
    \right)
\]
and    
\[
    D_{e_2} \text{ has type } 
    \left(
    \left\lceil\frac{p_2}{a_j}\right\rceil + 
    \left\lceil\frac{q_2}{a_i}\right\rceil,
    \eta(N_2) -
    \left\lceil\frac{p_2}{a_j}\right\rceil a_jd_i -  
    \left\lceil\frac{q_2}{a_i}\right\rceil a_id_j
    \right),
\]
and the intersection number $D_{e_1} \cdot D_{e_2}$ simplifies to the
left-hand side of \eqref{Eq:qsmoothintnum}.

Now suppose that conditions 1-3 hold.  We will show that the general
member is quasi-smooth.  It suffices to proceed according to irreducible
components present in the base locus of $|-K_{\wPP(\E)}|$.  If the base
locus is empty, then Proposition~\ref{P:qbertini} guarantees the existence
of a quasi-smooth member.  Condition 1 prevents there from being a fixed
component in the base locus.  Now suppose that the base locus has an
irreducible component $C$ of codimension 2 given by $\{X_{e_1} = X_{e_2}
=0 \}$.  Then we must check that it cannot be the case that all the
members fail to be quasi-smooth on $C$.  As described in the previous
paragraph, condition 3.(b) prevents this from happening.  Note that we
also need condition 1 to ensure that the two divisors $D_{e_1}$ and
$D_{e_2}$ from the previous paragraph intersect transversally. It
remains to check the case in which the base locus has an irreducible
component $C$ of codimension three.  Again, as described in the previous
paragraph, conditions 2.(b) and 2.(c) ensure that not every member fails
to be quasi-smooth on $C$.  This completes the proof.
\end{proof}

\begin{corollary}\label{C:qsmoothk3fib}
With the notation and assumptions of Theorem~\ref{T:qsmoothk3fib} the
general member of the linear system $|-K_{\wPP(\E)}|$ is quasi-smooth
and well-formed if and only if the following conditions hold.  
\begin{enumerate}
    \item For all $0 \leq i \leq 3$ there exists a monomial in $V_L$
    that does not involve $X_i$.  
    \item For all $0 \leq i \leq 3$
    \begin{enumerate}
        \item there exists a monomial of the form $X^p_i$ in $V_L$, or
        \item there exists an integer $j \neq i$ and a monomial $M \in
        V_L$ of the form $X_jX_i^{p}$ such that $\eta(M) = 0$, or
        \item there exist distinct integers $j_1 \neq i$ and $j_2 \neq
        i$ and two monomials $M_1$ and $M_2$ of the form
        $X_{j_1}X_{i}^{p}$ and $X_{j_2}X_{i}^{p}$ respectively such that
        $\eta(M_1) > 0$ and $\eta(M_2) > 0$.
    \end{enumerate}
    \item For all $0 \leq i < j \leq 3$
    \begin{enumerate}
        \item there exists a monomial of the form $X_i^pX_j^q \in V_L$
        or 
        \item there exist monomials of the form
        \[
        N_1 = X_i^{p_1}X_j^{q_1}X_{e_1} \in V_L
        \]
        and 
        \[
        N_2 = X_i^{p_2}X_j^{q_2}X_{e_2} \in V_L
        \]
        in $V_L$ such that $e_1$ and $e_2$ are distinct and
        such that
        \begin{multline}\label{Eq:qsmoothintnum1}
        \left(
        \left\lceil\frac{p_1}{a_j}\right\rceil +
        \left\lceil\frac{q_1}{a_i}\right\rceil
        \right)\eta(N_2) +
        \left(
        \left\lceil\frac{p_2}{a_j}\right\rceil +
        \left\lceil\frac{q_2}{a_i}\right\rceil
        \right)\eta(N_1) +\\
        \left(
        \left\lceil\frac{p_1}{a_j}\right\rceil 
        \left\lceil\frac{p_2}{a_j}\right\rceil +
        \left\lceil\frac{q_1}{a_i}\right\rceil
        \left\lceil\frac{q_2}{a_i}\right\rceil
        \right)(a_id_j - a_jd_i) = 0.
        \end{multline}
    \end{enumerate}
\end{enumerate}
\end{corollary}

\begin{proof}
This follows immediately from Theorem~\ref{T:qsmoothk3fib} and
Proposition~\ref{P:wp3qsmooth}.
\end{proof}

\section{Elliptically Fibered Calabi-Yau Threefolds}\label{S:ellfib}
Using techniques similar to those used in the proof of
Theorem~\ref{T:qsmoothk3fib} we can prove the analogous result for the
case of Calabi-Yau threefolds fibered over $\PP^2$ whose general fiber
is a genus one curve in $\PP(1,1,1)$, $\PP(1,1,2)$, or $\PP(1,2,3)$.
The statement of the full result follows.  It turns out that these
Calabi-Yau threefolds are nonsingular.  See \cite[Section~4.4]{jM06} for
details.

\begin{theorem}\label{T:smoothellfib}
Let $(d_0, d_1, d_2)$ be a sequence of integers that satisfy
\eqref{Eq:normform1} and \eqref{Eq:normform2}, and let $(a_0, a_1, a_2)$
be a sequence of positive integers satisfying \eqref{Eq:wform}.
Let $\E$ be the weighted locally free sheaf 
\[
    \E := \Oh_{\PP^2}(d_0) \oplus \Oh_{\PP^2}(d_1) \oplus
    \Oh_{\PP^2}(d_3)
\]
with weights $(a_0, a_1, a_2)$.  Let $L$ denote the linear system on
$\PP(a_0, a_1, a_2)$ obtained by restricting $|-K_{\wPP(\E)}|$ to a
fiber, and let $\eta$ be defined via \eqref{Eq:eta}  Then the linear
system $|-K_{\wPP(\E)}|$ has a well-formed nonsingular member if and
only if the following conditions hold.
\begin{enumerate}
    \item For all $i \in \{0,1,2\}$ there exists a monomial in $V_L$
    that does not involve $X_i$. 
    \item For all $i \in \{0,1,2\}$ either 
    \begin{enumerate}
        \item there exists a monomial of the form $X_i^{p}$ in $V_L$
        or 
        \item there exists an integer $j \neq i$ and a monomial $M$ of
        the form $X_i^{p}X_j$ such that $\eta(M) = 0$.
    \end{enumerate}
\end{enumerate}
\end{theorem}

There are 92 families of elliptically fibered Calabi-Yau threefolds over
$\PP^2$ arising as members of the anti-canonical linear system of
$\PP(a_0, a_1, a_2)$-bundle over $\PP^2$.  We list these families in the
appendix.

\section{Final remarks}\label{S:final}
It is natural to ask to what extent the Calabi-Yau threefolds we have
found are new.  From the point of view of mirror symmetry, they are not
new in the sense that they can be related to Batyrev's mirror symmetry
construction (see \cite{vB94}).  Although the varieties $\wPP(\E)$ are
not, in general, reflexive (i.e. Fano and Gorenstein), their Calabi-Yau
hypersurfaces have crepant resolutions that are anti-canonical
hypersurfaces in reflexive toric varieties, and these have all been
classified in \cite{mK00}.  The existence of these resolutions is
established in \cite{jM06} by checking that the Newton polyhedra of all
toric weighted projective bundles admitting quasi-smooth
anti-canonical hypersurfaces are reflexive.  

The varieties are interesting, however, for other reasons.  For example,
in \cite{mK02} fiber structures of toric Calabi-Yau hypersurfaces are
considered in the reflexive case.  It can be tricky to find fiber
structures from that point of view, and as far as the author is aware,
no such explicit description of weighted K3-fibered Calabi-Yau
threefolds as given here has appeared before, and it was not necessary
to have reflexive toric varieties to write them down.  Additionally, our
result gives a direct generalization of Reid's discovery of the ``famous
95'' families of weighted K3 hypersurfaces, and is thus interesting from
the point of view of explicit birational geometry.

\appendix
\section{Some data in the K3 case}\label{S:somedata}
We give here a partial list of K3-fibered Calabi-Yau hypersurfaces in
toric weighted projective bundles over $\PP^1$.  The list is in the
form of a hash table in \emph{Macaulay 2}.  (Note: The raw
\emph{Macaulay 2} output has been altered by inserting line-breaks.)
Each entry in the hash table consists of a \emph{key} and
an ASCII arrow ``\verb+=>+'' followed by a \emph{value}.  Each key is a
list of four integers, and each value is a list of lists of integers.
If \verb+{a,b,c,d}+ is a key and \verb+{e,f,g,h}+ is a member of the
list of values for \verb+{a,b,c,d}+, then the general member of the
linear system $|-K_{\wPP(\E)}|$ is quasi-smooth where $\E$ is the
weighted locally free sheaf
\[
    \Oh_{\PP^1}(e) \oplus \Oh_{\PP^1}(f) \oplus \Oh_{\PP^1}(g) \oplus
    \Oh_{\PP^1}(h)
\]
weighted by $(a,b,c,d)$.  The complete list of all 3,723 families and
details of the calculation appear in \cite{jM06}.

\footnotesize
\begin{verbatim}
{1, 1, 1, 1} => {{0, 0, 0, 0}, {0, 0, 0, 1}, {0, 0, 0, 2}, {0, 0, 1, 1}, 
                 {0, 0, 2, 2}, {0, 1, 1, 1}, {0, 1, 1, 2}, {0, 1, 1, 3}, 
                 {0, 1, 1, 4}}
{1, 1, 1, 2} => {{0, 0, 0, 0}, {0, 0, 0, 1}, {0, 0, 0, 2}, {0, 0, 1, 1}, 
                 {0, 0, 2, 2}, {0, 1, 1, 1}, {0, 1, 1, 2}, {0, 1, 1, 3}, 
                 {0, 1, 1, 4}, {0, 1, 1, 5}, {1, 1, 1, 1}, {0, 0, 1, 0}, 
                 {0, 0, 2, 0}, {0, 1, 2, 1}, {0, 1, 3, 1}, {1, 1, 2, 1}, 
                 {0, 1, 1, 0}, {0, 2, 2, 0}, {1, 2, 2, 1}, {1, 2, 3, 1}, 
                 {1, 2, 4, 1}, {1, 0, 0, 0}, {2, 0, 0, 0}, {1, 0, 1, 0}, 
                 {2, 0, 2, 0}, {1, 1, 1, 0}, {2, 1, 1, 0}, {3, 1, 1, 0}, 
                 {4, 1, 1, 0}, {5, 1, 1, 0}, {2, 1, 1, 1}, {2, 1, 2, 1}, 
                 {3, 1, 2, 1}}
{1, 1, 1, 3} => {{0, 0, 0, 0}, {0, 0, 0, 1}, {0, 0, 0, 2}, {0, 0, 1, 1}, 
                 {0, 0, 2, 2}, {0, 1, 1, 1}, {0, 1, 1, 2}, {0, 1, 1, 3}, 
                 {0, 1, 1, 4}, {0, 1, 1, 5}, {0, 1, 1, 6}, {1, 1, 1, 1}, 
                 {1, 1, 1, 2}, {1, 1, 2, 2}, {0, 0, 1, 0}, {0, 0, 2, 0}, 
                 {0, 1, 2, 1}, {0, 1, 1, 0}, {1, 0, 0, 0}, {2, 0, 0, 0}, 
                 {1, 0, 1, 0}}
{1, 1, 2, 2} => {{0, 0, 0, 0}, {0, 0, 0, 1}, {0, 0, 0, 2}, {0, 0, 1, 1}, 
                 {0, 0, 2, 2}, {0, 1, 1, 1}, {0, 1, 1, 2}, {0, 1, 1, 3}, 
                 {1, 1, 1, 1}, {1, 1, 1, 2}, {0, 1, 0, 1}, {0, 2, 0, 2}, 
                 {1, 2, 1, 2}, {1, 2, 1, 3}, {1, 2, 1, 4}, {1, 2, 1, 5}, 
                 {0, 1, 0, 0}, {0, 2, 0, 0}, {0, 2, 1, 1}, {0, 3, 1, 1}, 
                 {1, 2, 1, 1}, {1, 0, 0, 0}, {2, 0, 0, 0}, {1, 0, 0, 1}, 
                 {2, 0, 0, 2}, {1, 1, 0, 1}, {2, 1, 1, 1}, {2, 1, 1, 2}, 
                 {1, 1, 0, 0}, {0, 0, 1, 0}, {0, 0, 2, 0}, {1, 0, 1, 0}, 
                 {2, 0, 2, 0}, {1, 1, 1, 0}, {1, 1, 2, 0}, {1, 1, 3, 0}, 
                 {1, 1, 4, 0}, {1, 1, 5, 0}, {1, 1, 6, 0}, {1, 1, 2, 1}, 
                 {2, 1, 2, 1}, {2, 1, 3, 1}, {2, 1, 4, 1}}
{1, 1, 2, 3} => {{0, 0, 0, 0}, {0, 0, 0, 1}, {0, 0, 0, 2}, {0, 0, 1, 1}, 
                 {0, 0, 2, 2}, {0, 1, 1, 1}, {0, 1, 1, 2}, {0, 1, 1, 3}, 
                 {0, 1, 1, 4}, {1, 1, 1, 1}, {1, 1, 1, 2}, {1, 1, 1, 3}, 
                 {1, 1, 2, 2}, {1, 2, 2, 2}, {2, 2, 2, 2}, {0, 0, 1, 0}, 
                 {0, 0, 2, 0}, {0, 1, 2, 1}, {0, 1, 3, 1}, {1, 1, 2, 1}, 
                 {1, 1, 3, 2}, {1, 2, 3, 2}, {0, 1, 0, 1}, {0, 2, 0, 2}, 
                 {1, 2, 1, 2}, {1, 2, 1, 3}, {1, 2, 1, 4}, {1, 2, 1, 5}, 
                 {1, 2, 1, 6}, {1, 3, 1, 3}, {0, 1, 0, 0}, {0, 2, 0, 0}, 
                 {0, 2, 1, 1}, {1, 2, 1, 1}, {1, 4, 1, 3}, {2, 3, 2, 2}, 
                 {2, 4, 2, 2}, {0, 1, 1, 0}, {0, 2, 2, 0}, {1, 2, 2, 1}, 
                 {1, 3, 3, 2}, {1, 3, 2, 1}, {1, 4, 2, 1}, {1, 5, 2, 1}, 
                 {1, 4, 3, 2}, {1, 0, 0, 0}, {2, 0, 0, 0}, {1, 0, 0, 1}, 
                 {2, 0, 0, 2}, {1, 1, 0, 1}, {2, 1, 0, 1}, {3, 1, 0, 1}, 
                 {2, 1, 1, 1}, {2, 1, 1, 2}, {3, 1, 1, 3}, {4, 1, 1, 3}, 
                 {2, 2, 1, 2}, {3, 2, 1, 2}, {3, 2, 2, 2}, {4, 2, 2, 2}, 
                 {1, 1, 0, 0}, {2, 2, 0, 0}, {3, 2, 1, 1}, {1, 0, 1, 0}, 
                 {2, 0, 2, 0}, {1, 1, 1, 0}, {2, 1, 1, 0}, {3, 1, 1, 0}, 
                 {4, 1, 1, 0}, {2, 1, 2, 1}, {3, 1, 2, 1}, {4, 1, 2, 1}, 
                 {5, 1, 2, 1}, {1, 1, 0, 2}, {1, 1, 0, 3}, {1, 1, 0, 4}, 
                 {1, 1, 0, 5}, {1, 1, 0, 6}, {1, 1, 0, 7}, {2, 1, 1, 3}, 
                 {2, 1, 1, 4}, {2, 1, 1, 5}}
\end{verbatim}
\normalsize

\section{Data in the elliptic case}\label{S:elldata}
In the section we give the complete list of the elliptically fibered
Calabi-Yau threefolds discussed in Section~\ref{S:ellfib}.  This list is
to be interpreted in a manner similar to the list in
Appendix~\ref{S:somedata}.

\footnotesize
\begin{verbatim}
{1, 1, 1} => {{0, 0, 0}, {0, 0, 1}, {0, 0, 2}, {0, 0, 3}, {0, 1, 1}, 
              {0, 1, 2}, {0, 1, 3}, {0, 2, 3}}                
{1, 1, 2} => {{0, 0, 0}, {0, 0, 1}, {0, 0, 2}, {0, 0, 3}, {0, 1, 1}, 
              {0, 1, 2}, {0, 1, 3}, {0, 2, 3}, {0, 3, 3}, {0, 3, 4}, 
              {0, 3, 5}, {0, 3, 6}, {0, 3, 7}, {0, 3, 8}, {0, 3, 9}, 
              {0, 3, 10}, {0, 3, 11}, {0, 3, 12}, {1, 1, 1}, {0, 1, 0}, 
              {0, 2, 0}, {0, 3, 0}, {0, 2, 1}, {0, 3, 1}, {0, 3, 2}, 
              {0, 4, 3}, {0, 5, 3}, {0, 6, 3}, {1, 2, 1}, {1, 3, 1}, 
              {1, 1, 0}, {1, 2, 0}, {2, 2, 1}}
{1, 2, 3} => {{0, 0, 0}, {0, 0, 1}, {0, 0, 2}, {0, 0, 3}, {0, 1, 1}, 
              {0, 1, 2}, {0, 1, 3}, {0, 2, 3}, {0, 3, 3}, {0, 3, 4}, 
              {0, 3, 5}, {0, 3, 6}, {0, 3, 7}, {0, 3, 8}, {0, 3, 9}, 
              {1, 1, 1}, {1, 1, 2}, {1, 1, 3}, {1, 1, 4}, {1, 2, 2}, 
              {2, 2, 2}, {0, 1, 0}, {0, 2, 0}, {0, 3, 0}, {0, 2, 1}, 
              {0, 3, 1}, {0, 3, 2}, {0, 4, 3}, {0, 5, 3}, {0, 6, 3}, 
              {1, 2, 1}, {1, 3, 1}, {1, 3, 2}, {1, 4, 2}, {2, 3, 2}, 
              {1, 0, 1}, {1, 0, 2}, {2, 1, 2}, {2, 1, 3}, {1, 0, 0}, 
              {2, 0, 0}, {3, 0, 0}, {2, 0, 1}, {2, 1, 1}, {3, 1, 1}, 
              {3, 1, 2}, {3, 2, 2}, {1, 1, 0}, {1, 2, 0}, {2, 2, 1}, 
              {2, 1, 0}}
\end{verbatim}
\normalsize

\bibliographystyle{amsplain}
\bibliography{bib}
\end{document}